\documentclass{article}[12pt]
\usepackage{latexsym, amsmath, amsfonts, amssymb}
\usepackage{hyperref}
\usepackage{color}

\newcommand{\nc}{\newcommand}
\nc{\ga}{\gamma} \nc{\di}{\displaystyle}
\nc{\ek}{\protect\\[1ex]}
\nc{\N}{{\mathbb N}} \nc{\R}{{\mathbb R}} \nc{\Z}{{\mathbb Z}}
\nc{\La}{\Lambda} \nc{\la}{\lambda} \nc{\da}{\delta}
\nc{\Da}{\Delta} \nc{\na}{\nabla} \nc{\vp}{\varphi} \nc{\si}{\sigma}
\nc{\Si}{\Sigma} \nc{\al}{\alpha} \nc{\be}{\beta} \nc{\om}{\omega}
\nc{\Om}{\Omega} \nc{\pa}{\partial} \nc{\ti}{\times}
 \nc{\ve}{\varepsilon} \nc{\ra}{\rightarrow} \nc{\Ra}{\Rightarrow}
\nc{\ran}{\rangle} \nc{\lan}{\langle}

 \nc{\eq}[1]{\mbox{\rm{(\ref{E#1})}}}
\nc{\qed}{\mbox{}\nolinebreak\hfill \rule{2mm}{2mm}}
 \nc{\ha}{\frac{1}{2}}
\nc{\hra}{\hookrightarrow} \nc{\supp}{{\rm supp}\,}
\nc{\curl}{\text{curl}\,} \nc{\dense}{\hra^{\hspace{-3mm}d\,\,}}

\newtheorem{lem}{Lemma}[section]
\newtheorem{theo}[lem]{Theorem}

\newtheorem{rem}[lem]{Remark}

\renewcommand{\div}{{\rm{div}}\,}

\numberwithin{equation}{section} 

\title{Global regularity for the Navier-Stokes
 equations with application to global solvability\\ for the Euler equations}

\author{Myong-Hwan Ri\\
\small Institute of Mathematics, State Academy of Sciences, DPR
 Korea}
\date{}

  \allowdisplaybreaks[4]

\begin{document}
\bibliographystyle{alpha}

\maketitle%

\begin{abstract}
We show that any Leray-Hopf weak solution to the $d$-dimensional
Navier-Stokes equations $(d\geq 3)$ with initial values $u_0\in
H^{s}(\R^d)$, $s\geq -1+\frac{d}{2}$, belongs to $L^\infty(0,\infty;
H^{s}(\R^d))$ and thus it is globally regular. For the proof, first,
we construct a supercritical space which has very sparse inverse
logarithmic weight in the frequency domain, compared to the critical
homogeneous Sobolev $\dot{H}^{-1+d/2}$-norm.
 Then we obtain the energy estimates of high
frequency parts of the solution which involve the supercritical norm
as a factor of the upper bounds. Finally, we superpose the energy
norm of high frequency parts of the solution to get estimates of the
critical and subcritical norms independent of the viscosity
coefficient for the weak solution via the re-scaling argument.

As a direct application via the argument of vanishing viscosity, we
obtain global existence of a solution for the incompressible Euler
equations with initial values in $H^s(\R^d)$, $s\geq
-1+\frac{d}{2}$, which is unique when $s>1+\frac{d}{2}$.
\end{abstract}


\let\thefootnote\relax\footnote{\hspace{-0.3cm}
E-mail address\,$:$ math.inst@star-co.net.kp (Myong-Hwan Ri)}


\section{Introduction and main results}
Let us consider the Cauchy problem for the Navier-Stokes equations:
\begin{equation}
\label{E1.1}
\begin{array}{rl}
     u_t-\nu\Da u  + (u\cdot\na)u+\na p  = 0 \,\, &\text{in }(0,\infty)\ti\R^d,\ek
      \div u = 0 \,\, &\text{in }(0,\infty)\ti\R^d,\ek
u(0,x)=u_0\, &\text{in }\R^d.
\end{array}
\end{equation}

Since Leray \cite{Le34} proved existence of a global weak solution
(Leray-Hopf weak solution) $u$ to \eq{1.1} such that
$$
\label{E1.3} u\in L^2(0,\infty; H^1(\R^d))\cap L^\infty(0,\infty;
L^2(\R^d))
 $$
 satisfies \eq{1.1} in a weak sense and the {\it energy inequality}
\begin{equation}
\label{EEI}
\begin{array}{l}\di\frac{1}{2}\|u(t)\|_2^2+\nu\int_0^t\|\na
u(\tau)\|_2^2\,d\tau\leq \frac{1}{2}\|u_0\|_2^2,\; \forall t\in
(0,\infty),
 \end{array}
\end{equation}
the problem of global regularity of the weak solution for $d\geq 3$
has been a consistent key issue.

There is a great number of articles on scaling-invariant regularity
criteria for weak solutions to \eq{1.1}.
 We recall that the Navier-Stokes equations are
invariant by the scaling $u_\la(t,x)\equiv \la u(\la^2t,\la
x),\la>0,$ and
 a critical space for the equations is the space whose norm is invariant with respect to
the scaling
$$ u(x)\mapsto\la u(\la x), \; \la>0.$$
The following embedding holds
 between critical spaces for the  $d$-dimensional Navier-Stokes equations:
$$
\begin{array}{l}
\dot{H}^{-1+d/2}\hra   L^d\hra \dot{B}^{-1+d/p}_{p,\infty}\hra
  BMO^{-1}
  \hra  \dot{B}^{-1}_{\infty,\infty},\; d \leq p<\infty.
  \end{array}$$
 Eskauraiza, Seregin and Sv\'erak in \cite{ESS03}(2003) proved that a 3D Leray-Hopf weak solution
 $u$ to \eq{1.1} is regular in $(0,T]$ if
 $$u\in L^\infty(0,T; L^3(\R^3)),$$
  employing the backward uniqueness property of parabolic equations.
By developing a profile decomposition technique and using the method of ``critical
elements" developed in \cite{KeMe06}-\cite{KeMe10}, Gallagher, Koch and Planchon proved
for  a mild solution $u$ to \eq{1.1} with initial values in $L^3(\R^3)$ that a potential
 singularity at $t=T$ implies $\lim_{t\ra
 T-0}\|u(t)\|_3=\infty$ in \cite{GKB13}(2013), and extended the result
  from $L^3(\R^3)$ to wider critical Besov spaces
$\dot{B}^{-1+3/p}_{p,q}(\R^3)$, $3<p,q<\infty$
   in \cite{GKB16}(2016).
  In \cite{Ri23}(2023), the author proved that the above results can be extended
 to the largest critical space $\dot{B}^{-1}_{\infty,\infty}(\R^d)$ for the case of general $d\geq 3$, that is,
 a Leray-Hopf weak solution to \eq{1.1} satisfying
  $$u\in L^\infty(0,T; \dot{B}^{-1}_{\infty,\infty}(\R^d))$$
 is regular in $(0,T]$ by developing a superposition method of energy norms of high frequency parts
 for the weak solution.

Concerning supercritical regularity criteria, we mention recent
results \cite{BaPr21} and \cite{BaPr22} by Barker and Prange that a
3D Leray-Hopf weak solution $u$ satisfying the local energy
inequality is regular in $(0,T]$ if
$$\sup_{0<t<T}\int_{\R^3} \frac{|u(x,t)|^3}
{\Big(\log\log\log\big(\big(\log(e^{e^{3e^e}}+|u(x,t)|)\big)^{1/3}\big)\Big)^\theta}dx<\infty,\forall
\theta\in (0,1),$$ based on a Tao's result \cite{Tao19, Tao21} on
quantitative bound of the critical $L^3$-norm near a possible
blow-up epoch for a 3D weak solution to \eq{1.1}; we also mention
Pan \cite{Pan16} and Seregin \cite{Se22} where logarithmically
supercritical regularity criteria for axisymmetric suitable weak
solutions to \eq{1.1} are obtained.

\par\medskip
In this paper, we prove global regularity of the Leray-Hopf weak
solutions to \eq{1.1} for initial values in $H^{s}(\R^d)$, $d\geq
3$, $s\geq -1+d/2$. More precisely, we have:
\begin{theo}
\label{T1.1} {\rm Let $u$ be a Leray-Hopf weak solution to \eq{1.1}
with $u_0\in H^{s}(\R^d)$, $d\geq 3$, $s\geq -1+d/2$, and $\div
u_0=0$.
 Then,
$$u\in
L^\infty(0,\infty;H^{s}(\R^d))\cap L^2(0,\infty;H^{s+1}(\R^d))
$$
and the estimate
 \begin{equation}
 \label{E1.2n}
\|u(t)\|_{H^{s}(\R^d))}^2+\nu\int_0^t\|u(\tau)\|_{H^{s+1}(\R^d))}^2\,d\tau\leq
C(s)\|u_0\|_{H^s(\R^d)}^2, \;\forall t\in (0,\infty),
\end{equation}
holds true with a constant $C(s)>0$ depending only on $s$ and
independent of $\nu$ and $d$. In particular, $u$ is globally regular
and unique.
 }
\end{theo}
\begin{rem}
 {\rm
 	(i) Combining the result of Theorem \ref{T1.1} with already known regularity theory of the Navier-Stokes equations, cf. e.g. \cite{Am00}, one can conclude the following statement: 
 	
 	If $u_0\in H^s(\R^d)$, $s\geq -1+d/2$, $\div u_0=0$, then the problem \eq{1.1} has a unique global strong solution
 	 $$u\in C([0,\infty), H^s(\R^d))\cap L^2(0,\infty; H^{s+1}(\R^d))$$
 	 satisfying the estimate \eq{1.2n}.  
 	
(ii) The statement of Theorem \ref{T1.1} that  the estimate constants of
higher-order norms are independent of the viscosity is very
important. This fact yields, by vanishing viscosity, a global
existence of solutions for Euler equations of ideal incompressible
fluids, see Theorem \ref{T1.3}.
 }
\end{rem}

We use the following notations. The sets of all natural numbers and
all integers are denoted by $\N$ and $\Z$, respectively. The usual
$L^p$-norm for $1\leq p\leq \infty$ is denoted by $\|\cdot\|_p$.
  The Fourier
transform of a function $u$ is given by ${\cal F}u\equiv
\hat{u}=(2\pi)^{-d/2}\int_{\R^d}{u}(t, x)e^{ix\cdot\xi}\,dx$ and
${\cal F}^{-1}u\equiv \check{u}$. For $s\in\R$, $\dot{H}^s(\R^d)$
and $H^s(\R^d)$ stand for the homogeneous and inhomogeneous Sobolev
spaces with norms
$$\|u\|_{\dot{H}^s(\R^d)}\equiv
 \| |\xi|^s\hat{u}(\xi)\|_2\quad\text{and}\quad \|u\|_{H^s(\R^d)}\equiv
 \|(1+|\xi|^2)^{s/2}\hat{u}(\xi)\|_2,$$ respectively.
 We use the notation
\begin{equation} \label{E2.1}
\begin{array}{l}
 u^k(t,x):=(2\pi)^{-d/2}\int_{|\xi|\geq k}{\mathcal F}{u}(t,
\xi)e^{ix\cdot\xi}\,d\xi,\ek u_k:=u-u^k,\;
  u_{h,k}:=u^{h}-u^{k}\quad \text{for } 0\leq h<k.
  \end{array}
\end{equation}
We do not distinguish between the spaces of vectorial functions and
scalar functions.

Let us roughly explain the ideas for the proof of Theorem \ref{T1.1}
for the particular case $d=3$ and $s=1/2$; the idea is exactly the
same for general $d\geq 3$ and $s$. As is well-known, the
cancelation property $((v\cdot\na)u,u)=0$ for suitably smooth $u,v$
with $\div v=0$ is the key to obtain global existence of the
Leray-Hopf weak solutions to \eq{1.1}. However, this innermost
property was not used, so long as the author knows, in most previous
works for regularity. Indeed, as long as the convection term
$(u\cdot\na)u$ is treated as merely a quadratic nonlinear term,
whatever improvement of its estimate could be made, global
regularity of the weak solution is obtained under a smallness
condition of a critical norm of initial values or under additional
scaling-invariant conditions on the weak solution itself.

 In order to circumvent such situations, we use essentially the
 cancelation property $((u\cdot\na)u^k,u^k)=0$ by testing the momentum equation with high frequency parts $u^k$,
 $k\in\N$. Then, in the (local) time interval where the weak solution is regular we
 have
$$\frac{d}{2dt}\|u^k\|_2^2+\nu\|\na u^k\|_2^2\leq c\|u\|_{X}\|\na
  u^{k/2}\|_2^2,\;\forall k\in\N,$$
usually with a critical space $X$ or, through a more refined
observation,
$$\frac{d}{2dt}\|u^k\|_2^2+\nu\|\na u^k\|_2^2\leq c(k)\|u\|_{X_1}\|\na
  u^{k/2}\|_2^2,\;\forall k\in\N,$$ with a supercritical space $X_1$,
  where $c(k)$ depends on the topology of $X_1$.
Though this estimate, at the first glance, still seems of no special
use due to the factor $c(k)$, we pay attention to the fact that
$$\sum_{k\in\N}\|u^k\|_2^2\sim
\|u^1\|_{\dot{H}^{1/2}}^2,\quad \sum_{k\in\N}\|\na u^k\|_2^2\sim
\|\na u^1\|_{\dot{H}^{1/2}}^2$$ and
$$\sum_{k\in\N}c(k)\|\na  u^{k/2}\|_2^2 \lesssim \|\na u\|_2^2
+\sum_{n\in\N}\Big(\sum_{k=1}^n c(k)\Big)\|\na u_{n,n+1}\|_2^2.$$
Thus, if we can construct a supercritical space $X_1$ such that the
averaging condition $\sum_{k=1}^n c(k)\lesssim n$ is satisfied, then
we are led to
$$\frac{d}{2dt}\|u^1\|_{\dot{H}^{1/2}}^2+\nu\|\na u^1\|_{\dot{H}^{1/2}}^2\lesssim
 \|u\|_{X_1}(\|\na u\|_2^2+ \|\na u^{1}\|_{\dot{H}^{1/2}}^2),$$
yielding by energy inequality that
$$\ha\|u(t)\|_{\dot{H}^{1/2}}^2+\nu\int_0^t\|\na u\|_{\dot{H}^{1/2}}^2\,d\tau
 \leq \ha\|u_0\|_{H^{1/2}}^2
 +c\int_0^t \|u\|_{X_1}(\|\na u\|_2^2+ \|\na u\|_{\dot{H}^{1/2}}^2)d\tau$$
 for each fixed $t$. Hence, if additionally uniform smallness of $\|u(\tau)\|_{X_1}$ in $\tau
\in [0,t]$ can be guaranteed by a suitable re-scaling given for
fixed $t$, then in the last term of the right-hand side the part
$c\int_0^t \|u\|_{X_1} \|\na u\|_{\dot{H}^{1/2}}^2d\tau$ can be
absorbed into the left-hand side and the remaining part can be
estimated using the energy inequality. Thus
$\|u(t)\|_{\dot{H}^{1/2}}\lesssim \|u_0\|_{H^{1/2}}$ is obtained
irrespective of the viscosity and $t$. Remark here that the
re-scaling parameter should depend on $t$, but it does not matter.

We successfully construct a supercritical space $X_1$ with the
above-mentioned averaging condition and the uniform smallness
property by re-scaling, which has a very sparse inverse logarithmic
weight in the frequency domain compared to the critical
$\dot{H}^{-1+d/2}$-norm, Section 2.

The proof of Theorem \ref{T1.1} is given in Section 3.

The regularity estimate for the Navier-Stokes equations, see
\eq{1.2n}, which are uniform with respect to the viscosity
coefficient $\nu$ enables us to obtain easily global solvability for
the incompressible Euler equations, Section 4.

\section{A frequency-weighted scaling-variant space}

In this section, we introduce a supercritical space, the norm of
which is quite slightly weaker than that of
$\dot{H}^{-1+d/2}(\R^d)$.
\par\medskip
Let
$$\Da_j u= (\chi_j(\xi)\hat{u}(\xi))^\vee\;\text{for }j\in\Z, u\in {\cal
S}'(\R^d) \;\text{with }\hat{u}\in L_{loc}^1(\R_\xi^d),
$$ where $\chi_j$ is the characteristic function of the
set $\{\xi\in\R^d:  2^{j-1}\leq |\xi|< 2^{j}\}$ and ${\cal
S}'(\R^d)$ is the space of tempered distributions.

It immediately follows by Plancherel's theorem that for $s\in\R$ the
homogeneous Sobolev norm $\|u\|_{\dot{H}^s(\R^d)}$ is equivalent to
the norm
\begin{equation}
 \label{E3.1}
\big(\sum_{j\in\Z}\|2^{sj}\Da_j u\|^2_{L^2(\R^d)}\big)^{1/2}.
 \end{equation}

\medskip
 Let an infinite sequence
$\{a(j)\}_{j\in\Z}\subset \R$ be such that
 \begin{equation}
 \label{E2.2}
a(j):=\left\{\begin{array}{cl}  \log_2 j & \text{if
}j=i+2^{2^k}\quad\text{for some }k\in\N, i\in\Z, -k\leq i\leq k,\ek
                1      & \text{else}.
             \end{array}
      \right.
      \end{equation}
Let us define the space $X_1$ by
 \begin{equation}
 \label{E2.3n}
 \begin{array}{l}
X_1:=\{v\in \mathcal S'(\R^d): \hat{v}(\xi)\in
L^1_{loc}(\R_\xi^d),\;\{ 2^{j(-1+d/2)}a^{-1}(j)\|\Da_j
v\|_2\}_{j\in\Z}\in l_2\},\ek
 \|v\|_{X_1}:=\Big\|\{2^{j(-1+d/2)}a^{-1}(j)\|\Da_j v\|_2\}_{j\in\Z}\Big\|_{l_2},
 \end{array}
 \end{equation}
 where the sequence $\{a(j)\}$ is given by \eq{2.2}.

Obviously, $a(j)=1$ for $j\leq 0$ and $\dot{H}^{-1+d/2}(\R^d)$ is
densely embedded into $X_1$, and
$$\|v\|_{X_1} \leq \|v\|_{\dot{H}^{-1+d/2}(\R^d)},\forall v\in
\dot{H}^{-1+d/2}(\R^d).$$

We pursue supercritical properties of the space $X_1$ via the next
lemma.

\begin{lem}
 \label{L4.2n}
 {\rm
 (i) Let $v\in \dot{H}^{-1+d/2}(\R^d)$. Then,
$$\begin{array}{l}\forall \ve>0, \exists l_0:=\max\{M+1,
\log_2\big[\frac{\log_2M\|v\|_{X_1}}{\ve}\big]+1\}>0,\ek
  \|\la v(\la\cdot)\|_{X_1}\leq \ve,\quad  \forall
\la=2^{2^{2^l}} (l\geq l_0),
\end{array}$$
 where $M>0$ is such that
\begin{equation}
 \label{E4.4}
\sum_{|j|\geq M} 2^{2j(-1+d/2)} \|\Da_j v\|_2^2\leq \frac{\ve^2}{2}.
 \end{equation}

(ii) Let $u\in C\big([0,T],\dot{H}^{-1+d/2}(\R^d)\big)$ with
$0<T<\infty$. Then, for any sufficiently small $\ve>0$ there exists
$l_0>0$ such that
$$\begin{array}{l}
  \|\la u(t, \la\cdot)\|_{X_1}\leq \ve,\quad  \forall t\in [0,T],\;\forall
\la=2^{2^{2^l}} (l\geq l_0).
\end{array}$$

}
\end{lem}
 {\bf Proof:} -- {\it Proof of (i)}:
Let $v\in \dot{H}^{-1+d/2}(\R^d)$ and fix $\ve>0$ arbitrarily. Then,
in view of
$$ \|v\|_{\dot{H}^{-1+d/2}}\sim \Big\|\{2^{j(-1+d/2)}\|\Da_j v\|_2\}_{j\in\Z}\Big\|_{l_2},$$
  there is $M=M(\ve,v)\in \N$ satisfying \eq{4.4}.

 In the proof, for the moment, we use a short notation $v_\la:=\la v(\la\cdot)$.
  By definition of the space $X_1$ we have
 $$\|v_\la\|_{X_1}= \Big\|\{2^{j(-1+d/2)}a^{-1}(j)\|\Da_j v_\la\|_2\}_{j\in \Z}\Big\|_{l_2}.$$
  Observe that $({\mathcal F} v_\la)(\xi)= \la ({\mathcal F}
v(\la\cdot))(\xi)=\la^{1-d}(\mathcal Fv)(\frac{\xi}{\la})$ for all
$\la>0$. Then, for all $j\in\Z$
  $$\begin{array}{l}\Da_j v_\la (x)
   =\di\frac{1}{(2\pi)^{d/2}}\int_{2^{j-1}\leq |\xi|<2^j} {\cal F}v_\la(\xi)
e^{ix\cdot\xi}\,d\xi\\[3ex]
  \qquad =\di\frac{\la^{1-d}}{(2\pi)^{d/2}}\int_{2^{j-1}\leq |\xi|<2^j} {\cal
F}v(\frac{\xi}{\la}) e^{ix\cdot\xi}\,d\xi\\[3ex]
 \qquad =\di\frac{\la}{(2\pi)^{d/2}}\int_{\la^{-1}2^{j-1}\leq |\eta|<\la^{-1}2^j} {\cal
F}v(\eta) e^{i \la x\cdot\eta} \,d\eta.
  \end{array}$$
Hence, if $\la\equiv 2^{2^{2^l}}$,  $l\in\N$, then
 \begin{equation}
 \label{E2.6n}
 \begin{array}{l}\Da_j v_{\la}(x)=\di\frac{\la}{(2\pi)^{d/2}}
  \int_{2^{j-2^{2^l}-1}\leq |\eta|<2^{j-2^{2^l}}} {\cal
F}v(\eta) e^{i \la_0 x\cdot\eta}\,d\eta \\[3ex]
 \qquad =\la (\Da_{j-2^{2^l}}v)(\la x),
 \end{array}
 \end{equation}
 yielding
 $$\begin{array}{l}2^{j(-1+d/2)}\Da_j v_{\la}(x)=2^{j(-1+d/2)+2^{2^l}} (\Da_{j-2^{2^l}}v)(\la x),\;\forall j\in\Z.
 \end{array}$$
Note that
 $$\|(\Da_{j-2^{2^l}}v)(\la \cdot)\|_2=\la^{-d/2}\|\Da_{j-2^{2^l}}v\|_2.$$
Therefore, by \eq{2.6n} we have
 \begin{equation}
 \begin{array}{rcl}
\di 2^{j(-1+d/2)}\|\Da_j v_{\la}\|_2&=& 2^{j(-1+d/2)+2^{2^l}}
\|(\Da_{j-2^{2^l}}v)(\la \cdot)\|_2\ek
 & = &\di 2^{(-1+d/2)(j-2^{2^l})}\|\Da_{j-2^{2^l}}v\|_2,\;\forall j\in\Z,
 \end{array}
 \end{equation}
and, in view of \eq{4.4},
$$\begin{array}{l}
 \|v_{\la}\|_{X_1}^2=\di\Big\|\big\{2^{j(-1+d/2)}a^{-1}(j)\|\Da_j v_{\la}\|_2\big\}_{j\in\Z}\Big\|_{l_2}^2\\[2ex]
 =\di\Big\|\big\{2^{(-1+d/2)(j-2^{2^l})}a^{-1}(j)
 \|\Da_{j-2^{2^l}}v\|_2\big\}_{j\in\Z}\Big\|_{l_2}^2\\[2ex]
 =\di\Big\|\big\{\frac{1}{a(j+2^{2^l})}2^{j(-1+d/2)}
 \|\Da_{j}v\|_2\big\}_{j\in\Z}\Big\|_{l_2}^2\\[2ex]
 \leq \di\Big\|\big\{\frac{1}{a(j+2^{2^l})}2^{j(-1+d/2)}
 \|\Da_{j}v\|_2\big\}_{|j|\leq M}\Big\|_{l_2}^2+\Big\|\big\{2^{j(-1+d/2)}
 \|\Da_{j}v\|_2\big\}_{|j|> M}\Big\|_{l_2}^2\\[2ex]
  \di\leq \frac{\ve^2}{2}+\di\Big\|\big\{\frac{a(j)}{a(j+2^{2^l})}2^{j(-1+d/2)}a^{-1}(j)
 \|\Da_{j}v\|_2\big\}_{|j|\leq M}\Big\|_{l_2}^2.
 \end{array}$$
 Here, for $l\geq M+1$, by definition of $a(j)$ we have
 $$\begin{array}{l}
  \di\big\|\{\frac{a(j)}{a(j+2^{2^l})}2^{j(-1+d/2)}a^{-1}(j)
 \|\Da_{j}v\|_2\}_{|j|\leq M}\big\|_{l_2}\ek\qquad
=  \di\big\|\{\frac{a(j)}{\log_2(j+2^{2^l})}2^{j(-1+d/2)}a^{-1}(j)
 \|\Da_{j}v\|_2\}_{|j|\leq M}\big\|_{l_2}\ek\qquad
 \leq \di\sup_{|j|\leq M}\frac{\log_2M}{\log_2(j+2^{2^l})}\cdot
 \big\|\{2^{j(-1+d/2)}a^{-1}(j)\|\Da_{j}v\|_2\}_{|j|\leq M}\big\|_{l_2}
 \\[3ex]
 \qquad \leq \di\frac{\log_2M}{2^{l-1}}\|v\|_{X_1}.
  \end{array}$$
Therefore, for all
 \begin{equation}
 \label{E2.11}
l\geq l_0:=\max\{M+1,
\log_2\big[\frac{\log_2M\|v\|_{X_1}}{\ve}\big]+1\}, \quad
\la=2^{2^{2^l}},
 \end{equation}
  we have
$$ \|v_{\la}\|_{X_1}^2\leq \frac{\ve^2}{2}+\di\Big(\frac{\log_2M}{2^{l-1}}\|v\|_{X_1}\Big)^2
  \leq \ve^2.
 $$

\par\bigskip
{\it -- Proof of (ii)}: Let $u\in C([0,T];\dot{H}^{-1+d/2}(\R^d))$
and put
$$f(t,s):=\|u_{1/s}(t)\|_{\dot{H}^{-1+d/2}}+\|u^s(t)\|_{\dot{H}^{-1+d/2}},\quad (t,s)\in [0,T]\ti [2,\infty).$$
Then,  $u\in C([0,T];\dot{H}^{-1+d/2}(\R^d))$ implies that $f\in
BC([0,T]\ti [2,\infty))$ and for all $t\in [0,T]$ the mapping $s \to
f(t,s)$ is nonincreasing in $s$ and $f(t,s)\rightarrow 0$ as $s\ra
\infty$. Therefore, for any sufficiently small $\ve >0$ and $t\in
[0,T]$ there exists $s=s(t)$ satisfying
$f(t,s(t))=\frac{\ve}{2^{(d-1)/2}}$, where
 one may choose  $s(t)$ such that
 \begin{equation}\label{E2.12}
\exists N>0,\; \forall t\in [0,T]:  \quad |s(t)|\leq N
 \end{equation}
  holds. In
fact, if \eq{2.12} does not hold true,  there is a sequence
$\{t_n\}\subset [0,T]$ satisfying $f(t_n,
s(t_n))=\frac{\ve}{2^{(d-1)/2}}$ and $s(t_n)\ra\infty$ as
$n\ra\infty$. Without loss of generality, we may assume $t_n\ra
\tilde{t}$ for some $\tilde{t}\in [0,T]$ for $n\ra \infty$. Then, in
view of $u\in C([0,T],\dot{H}^{-1+d/2}(\R^d))$, we are led to
$$\begin{array}{rl}
 \di\frac{\ve}{2^{(d-1)/2}}&=f(t_n, s(t_n))\leq |f(t_n, s(t_n))-f(\tilde{t}, s(t_n))|+f(\tilde{t},
 s(t_n))\ek
 &\leq \big|\|[u(t_n)]_{1/s(t_n)}\|_{{\dot{H}^{-1+d/2}}}-\|[u(\tilde{t})]_{1/s(t_n)}\|_{\dot{H}^{-1+d/2}}\big|\ek
&\qquad
+\big|\|[u(t_n)]^{s(t_n)}\|_{{\dot{H}^{-1+d/2}}}-\|[u(\tilde{t})]^{s(t_n)}\|_{\dot{H}^{-1+d/2}}\big|+f(\tilde{t},
 s(t_n))\ek
 &\leq \|[u(t_n)-u(\tilde{t})]_{1/s(t_n)}\|_{\dot{H}^{-1+d/2}}
 +\|[u(t_n)-u(\tilde{t})]^{s(t_n)}\|_{\dot{H}^{-1+d/2}}+f(\tilde{t},
 s(t_n))\ek
 &\leq \sqrt{2}\|u(t_n)-u(\tilde{t})\|_{\dot{H}^{-1+d/2}}+f(\tilde{t},
 s(t_n))\ek
 &\ra 0 \quad (n\ra\infty),
  \end{array}$$
which is a contradiction.

Consequently, if we put
$$M:=[\log_2N]+2\quad\text{(with $N$ in \eq{2.12})}$$
then, in view of $N\leq 2^{M-1}$, we have
$$\begin{array}{rcl}
 \di\sum_{|j|\geq M} 2^{2j(-1+d/2)} \|\Da_j u(t)\|_2^2 &=& 2^{d-2}\di\sum_{|j|\geq M} 2^{2(j-1)(-1+d/2)} \|\Da_j
 u(t)\|_2^2\ek
&\leq &
2^{d-2}\big(\|[u(t)]_{2^{-M+1}}\|_{\dot{H}^{-1+d/2}}^2+\|[u(t)]^{2^{M-1}}\|_{\dot{H}^{-1+d/2}}^2\big)\ek
&\leq &
2^{d-2}\big(\|[u(t)]_{1/N}\|_{\dot{H}^{-1+d/2}}^2+\|[u(t)]^{N}\|_{\dot{H}^{-1+d/2}}^2\big)\ek
&\leq &
2^{d-2}\big(\|[u(t)]_{1/N}\|_{\dot{H}^{-1+d/2}}+\|[u(t)]^{N}\|_{\dot{H}^{-1+d/2}}\big)^2\ek
&\leq & 2^{d-2}2^{1-d}\ve^2= \di\frac{\ve^2}{2}
 \end{array}$$
and thus \eq{4.4} holds true for $u(t)$ uniformly for all $ t\in
[0,T]$.

 Thus,  for
$$l\geq l_0:=\max\{M+1,
\log_2\big[\frac{\log_2M\|u\|_{C([0,T],
\dot{H}^{-1+d/2})}}{\ve}\big]+1\},$$ see \eq{2.11},  we have
$$ \|(u(t))_{\la}\|_{X_1}^2\leq \frac{\ve^2}{2}+\di\Big(\frac{\log_2M}{2^{l-1}}\|u(t)\|_{X_1}\Big)^2
  \leq \ve^2,\quad \forall t\in [0,T],
 $$
by the above proved fact (i) and $\|u(t)\|_{X_1}\leq
\|u(t)\|_{\dot{H}^{-1+d/2}}$.

 \hfill\qed

\section{Proof of Theorem \ref{T1.1}}

\begin{lem}
 \label{L4.1}
 {\rm
 Let
   \begin{equation}
  \label{E3.2}
b(j) \equiv 2^{-j-1}\sum_{i=1}^{j} 2^{i}a(i),\;j\in\N.\end{equation}
   Then, for all $j\in\N$
$$b(j)\leq \left\{\begin{array}{cl}  \log_2 j
     & \text{if }-k+2^{2^k}\leq j\leq 2^{2^k}+2k\quad\text{for some }k\in\N,\ek
                2      & \text{else}.
             \end{array}
      \right.$$
      }
\end{lem}
 {\bf Proof:}
 Note that $a(j)\geq 1$ and
$b(j+1)=\frac{1}{2}(b(j)+a(j+1))$ for all $j\in\N$.

 First we prove the lemma for $j=1\sim  2^{2^1}=4$. We have
 \begin{equation}
 \label{E2.13}
\begin{array}{l}
  b(1)=2^{-2}\cdot 2\cdot a(1)=\frac{1}{2}<a(1)=1,\ek
  b(2)=2^{-3}\cdot (2+2^2)=\frac{3}{4}<a(2),\ek
  b(3)=\frac{b(2)+a(3)}{2}=\frac{3+4\log_2 3}{8}\leq \log_2 3=a(3),\ek
b(4)=\frac{b(3)+a(4)}{2}\leq 2= a(4).
 \end{array}
 \end{equation}

Next, let us fix any $k\in\N$ and prove the lemma for $2^{2^k}<
j\leq 2^{2^{k+1}}$. Let $j_1=2^{2^k}$ and assume in view of
\eq{2.13} that $b(j_1)\leq a(j_1)$ which, at the moment, is
satisfied for $k=1$.
 Then, for
all $l=1,2,\ldots,k$ we have
$$\begin{array}{l}
  b(j_1+1)=\frac{b(j_1)}{2}+\frac{a(j_1+1)}{2}\leq \frac{a(j_1)}{2}+\frac{a(j_1+1)}{2}\leq a(j_1+1)=\log_2(j_1+1),\ek
   b(j_1+2)=\frac{b(j_1+1)}{2}+\frac{\cdot a(j_1+2)}{2}\leq \frac{a(j+1)}{2}+\frac{a(j_1+2)}{2}
    \leq a(j_1+2)=\log_2(j_1+2),\ek
     \hspace{3cm}\cdots\quad\cdots\quad\cdots \ek
  b(j_1+l)=\frac{b(j_1+l-1)}{2}+\frac{a(j_1+l)}{2}\leq a(j_1+l)=\log_2(j_1+l).
  \end{array}$$
Moreover, for $l\in\N $ with $1\leq l \leq 2^{2^{k+1}}-(2k+2^{2^k})$
we have
$$\begin{array}{l}
  b(j_1+k+1)=\frac{b(k+2^{2^k})}{2}+\frac{a(k+2^{2^k}+1)}{2}\leq \frac{a(k+2^{2^k})}{2}+\frac{1}{2}
   \leq \frac{\log_2(k+2^{2^k})}{2}+\frac{1}{2},\ek
  b(j_1+k+2)=\frac{b(k+2^{2^k}+1)}{2}+\frac{a(k+2^{2^k}+2)}{2}\leq \frac{\log_2(k+2^{2^k})}{4}
       +\frac{1}{4}+\frac{1}{2},\ek
     \hspace{3cm}\cdots\quad\cdots\quad\cdots \ek
  b(j_1+k+l)=\frac{b(k+2^{2^k}+l-1)}{2}+\frac{a(k+2^{2^k}+l)}{2}\leq
      \frac{\log_2(k+2^{2^k})}{2^{l}}+\frac{1}{2^{l}}+\cdots\frac{1}{4}+\frac{1}{2}.
  \end{array}$$
Hence, if $1\leq l <\log_2\log_2(k+2^{2^k})$, then we have
 $$\begin{array}{l} b(j_1+k+l)\leq \log_2(k+2^{2^k})
  \end{array}$$
  and, if $\log_2\log_2(k+2^{2^k})\leq l \leq 2^{2^{k+1}}-(2k+2^{2^k})$, then we have
 $$\begin{array}{l} b(j_1+k+l)\leq 2.
  \end{array}$$
Note that $k<\log_2\log_2(k+2^{2^k})<k+1$ for any $k\in\N$.
Consequently,
 \begin{equation}\label{E2.10}
\begin{array}{l}
 b(j_1+k+l)\leq \log_2(k+2^{2^k}), \quad \text{for }l=1,\ldots, k,\ek
b(j_1+k+l)\leq 2, \quad\text{for }l=k+1,\ldots,
2^{2^{k+1}}-(2k+2^{2^k}).
  \end{array}
  \end{equation}

Finally, we shall prove
 \begin{equation}
 \label{E2.14}
b(j)\leq \log_2 j\quad\text{for }2^{2^{k+1}}-k+1\leq j\leq
2^{2^{k+1}}.
 \end{equation}
In fact, we have
$$\begin{array}{l}
b(2^{2^{k+1}}-k+1)=\frac{b(2^{2^{k+1}}-k)}{2}+\frac{a(2^{2^{k+1}}-k+1)}{2}\ek
  \qquad\qquad \leq 1+\frac{\log_2(2^{2^{k+1}}-k+1)}{2}\leq \log_2(2^{2^{k+1}}-k+1), \ek
 b(2^{2^{k+1}}-k+2)=\frac{b(2^{2^{k+1}}-k+1)}{2}+\frac{a(2^{2^{k+1}}-k+2)}{2}\ek
 \qquad \qquad \leq
  \frac{\log_2(2^{2^{k+1}}-k+1)}{2}+\frac{\log_2(2^{2^{k+1}}-k+2)}{2}
   \leq \log_2(2^{2^{k+1}}-k+2),\ek
 \hspace{3cm}\cdots\quad\cdots\quad\cdots \ek
 b(2^{2^{k+1}})=\frac{b(2^{2^{k+1}}-1)}{2}+\frac{a(2^{2^{k+1}})}{2}\ek
 \qquad \qquad \leq
  \frac{\log_2(2^{2^{k+1}}-1)}{2}+\frac{\log_22^{2^{k+1}}}{2}
   \leq 2^{k+1}=\log_2{2^{2^{k+1}}}=a(2^{2^{k+1}}).
\end{array}$$
Here, in particular, $b(2^{2^{k+1}})\leq a(2^{2^{k+1}})$ together
with $b(4)\leq a(4)$ implies by induction that the assumption
$b(j_1)\leq a(j_1)$ is satisfied for all $k\in\N$.

Thus, by \eq{2.13} -- \eq{2.14}, the proof of the lemma is complete.

 \hfill\qed

In the below we use the notation
 \begin{equation}
 \label{E4.1}
 S(n)=\{l\in\N:\; l\leq n, -k+2^{2^k}\leq l\leq 2^{2^k}+2k,\exists k\in\N\},\;n\in\N.
  \end{equation}
\begin{lem}
 \label{L4.2}
 {\rm
The cardinal number $|S(n)|$ of the set $S(n)$ is estimated by
 $$|S(n)|\leq \di\frac{(3\log_2\log_2n +5)\log_2\log_2n}{2} .$$
 }
\end{lem}
 {\bf Proof:} Suppose that  $2^{2^{k-1}}\leq n <2^{2^k}$ for some $k\in\N$.
  Then, by definition of $S(n)$, we have
$$\begin{array}{rcl}
 |S(n)| & \leq  & \di\sum_{l=1}^{k-1} (3l+1)=\frac{(3k+2)(k-1)}{2}
 \\[3ex]
 &< & \di\frac{(3\log_2\log_2n +5)\log_2\log_2n}{2}.
\end{array}
$$
Thus, the lemma is proved. \hfill\qed

\par\bigskip\noindent
 {\bf Proof of Theorem \ref{T1.1}:}
Since $u_0\in H^{-1+d/2}(\R^d)$, the problem \eq{1.1} has a unique
local strong (mild) solution  $u\in C([0,T),
\dot{H}^{-1+d/2}(\R^d))$ (This fact follows by Kato's approach
 just changing $L^d(\R^d)$ in \cite{Ka84} with $\dot{H}^{-1+d/2}(\R^d))$,
cf. \cite{KaFu62}) and it coincides in its existence interval with
the corresponding Leray-Hopf weak solution (\cite{Am00}, Theorem
0.8). Moreover, it is smooth in $(0,T)$, as is well known (cf. e.g.
\cite{Ga00}). Let $(0,T)$ be the maximal interval where a Leray-Hopf
weak solution $u$ to \eq{1.1} with $u_0\in H^{-1+d/2}(\R^3)$ is
smooth. It is well-known that energy inequality \eq{EI} (even energy
equation) is satisfied in $(0,T)$.

 The remaining proof is divided into two
 steps.

\par\medskip
{\bf Step 1. Energy estimate of high frequency parts:}

Recall the notation \eq{2.1}. By $L^2$-scalar product to the
momentum equation of \eq{1.1} with $u^k(t)$, $k\geq 1$, $t\in
(0,T)$, we have
 \begin{equation}
 \label{E2.18}
\begin{array}{l}
\frac{d}{2dt}\|u^k\|_2^2+\nu\|\na u^k\|_2^2=-((u\cdot\na)u_k,
u^k)\ek
 \quad =-((u_k\cdot\na)u_k+(u^k\cdot\na)u_k, u^k)\ek
 \quad
=-((u_k\cdot\na)u_{k/2,k},u^{k})-((u_{k/2,k}\cdot\na)u_{k/2},u^{k})-((u^k\cdot\na)u_k,u^k)
 \;\text{in } (0,T)
\end{array}
 \end{equation}
  in view of  $u_k=u_{k/2}+u_{k/2,k}$ and
$$ \begin{array}{l}
 (u_k, u^l)_2=0,\quad \supp{\widehat{(u_ku_l)}}\subset \{\xi\in\R^d: |\xi|\leq k+l\},\;\quad k\leq l<\infty.
 \end{array}
 $$

 In the right-hand side of \eq{2.18} we have
$$\begin{array}{l}
|((u_{k/2,k}\cdot\na)u_{k/2},u^{k})+(u^k\cdot\na)u_k,u^k)|
 \leq  (\|\na u_{k/2}\|_\infty+\|\na u_k\|_\infty)\|u^{k/2}\|_2^2
  \end{array}$$
  by H\"older's inequality.
 On the other hand, we have
$$ \begin{array}{rcl} |((u_k\cdot\na)u_{k/2,k},u^{k})|
 &\leq &  \|u_k\|_\infty\|\na u_{k/2,k}\|_2\|u^k\|_2\ek
  &\leq &  k\|u_k\|_\infty\|u^{k/2}\|_2^2.
 \end{array}
 $$
 Therefore, we get from \eq{2.18} that
 \begin{equation}
 \label{E2.19}
 \begin{array}{l}
 \frac{d}{2dt}\|u^k\|_2^2+\nu\|\na u^k\|_2^2
      \leq (\|\na u_{k/2}\|_\infty+\|\na u_k\|_\infty+k\|u_k\|_\infty)\|u^{k/2}\|_2^2
      \quad\text{in }(0,T).
      \end{array}
 \end{equation}

Let
 \begin{equation}
 \label{E4.3}
j_0\equiv j_0(k)=\lceil\log_2 k\rceil+1.
 \end{equation}
Note that
 $$\begin{array}{rcl} \|\widehat{\Da_ju_k}\|_1
 &\leq & \|\widehat{\Da_ju_k}\|_2\cdot |\{2^{j-1}\leq |\xi|\leq
 2^j\}|^{1/2}\ek
  &\leq &  c_02^{dj/2}\|\Da_ju_k\|_2\ek
   &\leq &  c_02^{dj/2}\|\Da_ju\|_2,\;\forall j\in\Z,
 \end{array}$$
 where $c_0=\si^{1/2}$ and $\si$ is the volume of the unit ball of $\R^d$. Hence, we get by
  the definition of the space $X_1$
 that
 $$\begin{array}{rl}
 \|u_k\|_\infty&\leq (2\pi)^{-d/2}\|\widehat{u_k}\|_1
  =(2\pi)^{-d/2}\big\|\di\sum_{j=-\infty}^{j_0}\widehat{\Da_j u_k}\big\|_1\ek
  &\leq (2\pi)^{-d/2}
  \di\sum_{j=-\infty}^{j_0}\|\widehat{\Da_ju_k}\|_1\ek
  &\leq c_0(2\pi)^{-d/2}\di\sum_{j=-\infty}^{j_0} 2^{j}a(j)
     (a^{-1}(j)2^{j(-1+d/2)}\| \Da_ju\|_2)\ek
 &\leq
 \di c_0(2\pi)^{-d/2}\big(\sum_{j=-\infty}^{j_0} 2^{j}a(j)\big)
   \cdot \sup_{j\leq j_0}\{a^{-1}(j)2^{j(-1+d/2)}\big\| \Da_ju \big\|_2\}\ek
 &\leq c_0(2\pi)^{-d/2}(2+  b(j_0(k)) 2^{j_0+1})  \big\|u\big\|_{X_1}\ek
   &\leq 8c_0(2\pi)^{-d/2}b(j_0(k))k\|u\|_{X_1}
  \end{array}$$
  due to $k\geq 1$ and $b(j)\geq 1/2$ for all $j\in\N$, see \eq{2.2} and \eq{3.2}.
In the same way, we have
 $$\begin{array}{rl}
 \|\na u_k\|_\infty &\leq
 \di c_0(2\pi)^{-d/2}\big(\sum_{j=-\infty}^{j_0} 2^{j}a(j)\big)
   \cdot \sup_{j\leq j_0}\{a^{-1}(j)2^{j(-1+d/2)}\big\| \Da_j\na u \big\|_2\}\ek
 &\leq \di c_0(2\pi)^{-d/2}(2+  b(j_0(k)) 2^{j_0+1})\cdot 2^{j_0}
   \sup_{j\leq j_0}\{a^{-1}(j)2^{j(-1+d/2)}\big\| \Da_ju \big\|_2\}\ek
   &\leq 8c_0(2\pi)^{-d/2}b(j_0(k))k^2\|u\|_{X_1}
  \end{array}$$
and
 $$\begin{array}{rl}
 \|\na u_{k/2}\|_\infty &\leq
 \di c_0(2\pi)^{-d/2}\big(\sum_{j=-\infty}^{j_0-1} 2^{j}a(j)\big)
   \cdot \sup_{j\leq j_0-1}\{a^{-1}(j)2^{j(-1+d/2)}\big\| \Da_j\na u \big\|_2\}\ek
&\leq
 \di c_0(2\pi)^{-d/2}\big(\sum_{j=-\infty}^{j_0} 2^{j}a(j)\big)
   \cdot 2^{j_0-1}\sup_{j\leq j_0-1}\{a^{-1}(j)2^{j(-1+d/2)}\big\| \Da_j u \big\|_2\}\ek
   &\leq 4c_0(2\pi)^{-d/2}b(j_0(k))k^2\|u\|_{X_1}.
  \end{array}$$

 Thus,  we have by \eq{2.19} that
 \begin{equation}
 \label{E2.15}
\begin{array}{l}
\di  \frac{d}{2dt}\|u^{k}\|_2^2+\nu \|\na u^k\|_2^2
   \leq  C_1b(j_0(k))k^2\|u\|_{X_1}\|u^{k/2}\|_2^2\ek
\qquad\leq  4C_1b(j_0(k))\|u\|_{X_1}\|\na u^{k/2}\|_2^2, \,\forall
k\geq 1,\quad\text{in }(0,T),
\end{array}
 \end{equation}
with a generic constant $C_1>0$ depending on $d$ and independent of
$\nu$.

\par\bigskip
{\bf Step 2. Estimates of higher order norms for the solution:}

Adding up \eq{2.15} multiplied by $k^{s}$, $s\geq d-3,$ over all
$k\in \N$, we have
 \begin{equation}
 \label{E2.5n}
  \begin{array}{l}
 \di\frac{d}{2dt}\sum_{k\in\N}k^s\|u^{k}\|_2^2+\nu \sum_{k\in\N}
 k^s\|\na u^{k}\|_2^2  \leq  4C_1\|u(t)\|_{X_1}  \sum_{k\in\N} b(j_0(k))k^s\|\na u^{k/2}\|_2^2
   \end{array}
 \end{equation}
 in $(0,T)$.
With the notation $2\N:=\{2n: n\in\N\}$ we have
 \begin{equation}
 \label{E2.26nn}
\begin{array}{l}
  \di\sum_{k\in\N} b(j_0(k))k^s\|\na u^{k/2}\|_2^2\\[3ex]
   \di = \ha\|\na u^{1/2}\|_2^2+
   \sum_{k\in 2\N} \big(b(j_0(k))k^s\|\na u^{k/2}\|_2^2+b(j_0(k+1))(k+1)^s\|\na u^{(k+1)/2}\|_2^2\big)
\\[3ex]
   \di \leq  \|\na u^{1/2}\|_2^2+
   \sum_{k\in 2\N} \big(b(j_0(k))k^s+b(j_0(k+1))(k+1)^s\big)\|\na u^{k/2}\|_2^2
\\[3ex]
   \di =  \| \na u^{1/2}\|_2^2+
   \sum_{k\in \N} (b(j_0(2k))(2k)^s+b(j_0(2k+1)))(2k+1)^s\| \na u^{k}\|_2^2
\\[3ex]
\di =\|\na u^{1/2}\|_2^2+\sum_{k\in\N}
\big(b(j_0(2k))(2k)^s+b(j_0(2k+1))(2k+1)^s\big)\sum_{n\geq k}\|\na u_{n,n+1}\|_2^2\\[3ex]
\di =\|\na u^{1/2}\|_2^2+\sum_{n\in\N} \|\na u_{n,n+1}\|_2^2
\sum_{k=1}^n \big(b(j_0(2k))(2k)^s+b(j_0(2k+1))(2k+1)^s\big)
\\[3ex]
\di \leq \|\na u^{1/2}\|_2^2+\sum_{n\in\N} (3n)^s\|\na
u_{n,n+1}\|_2^2 \sum_{k=1}^n \big(b(j_0(2k))+b(j_0(2k+1))\big).
   \end{array}
  \end{equation}

Note that $j_0(2k)=\lceil \log_2k\rceil +2$, see \eq{4.3}, and
 $$\begin{array}{l}
\di\sum_{k=1}^n b(j_0(2k))\\[3ex]
 \quad= \di\sum_{k\leq n,j_0(2k)\notin S(\lceil\log_2n\rceil+2)} b(j_0(2k))
 +\sum_{k\leq n,j_0(2k)\in S(\lceil\log_2n\rceil+2)} b(j_0(2k)),
 \end{array}$$
where, by Lemma \ref{L4.1},
 $$\begin{array}{l}
\di\sum_{k\leq n,j_0(2k)\notin S(\lceil\log_2n\rceil+2)} b(j_0(2k))
\leq 2n.
 \end{array}$$
Moreover,  we have by \eq{4.1} and Lemma \ref{L4.2} that
$$\begin{array}{l}
\di\sum_{k\leq n,j_0(2k)\in S(\lceil\log_2n\rceil+2)} b(j_0(2k))\\[3ex]
 \qquad \qquad \di\leq  |S(\lceil\log_2n\rceil+2)|\cdot
           \max_{k\leq n, j_0(2k)\in S(\lceil\log_2n\rceil+2)} b(j_0(2k))\\[3ex]
   \qquad \qquad  \di\leq  |S(\lceil\log_2n\rceil+2)|\cdot
   \log_2(\lceil\log_2n\rceil+2)\\[3ex]
\qquad \qquad  \di\leq
\frac{3\log_2\log_2(\lceil\log_2n\rceil+2)+5}{2}\cdot
\log_2\log_2(\lceil\log_2n\rceil+2)\cdot\log_2(\lceil\log_2n\rceil+2)\\[3ex]
\qquad \qquad  \di \leq n,\quad \forall n\geq n_0,
 \end{array}$$
 for some generic number $n_0\in\N$.
Therefore, for all $n\geq n_0$ we have
 \begin{equation}
 \label{E3.6}
\begin{array}{l}
\di\sum_{k=1}^n b(j_0(2k)) \leq 2n+n\leq 3n.
 \end{array}
 \end{equation}
In the same way, we have
 \begin{equation}
 \label{E3.6n}
\begin{array}{l}
\di\sum_{k=1}^n b(j_0(2k+1)) \leq 3n,\quad \forall n\in\N, n\geq
n_0,
 \end{array}
 \end{equation}
with a generic number $n_0\in\N$  possibly larger than the former
$n_0$.


 Thus we get by
\eq{2.26nn}-\eq{3.6n} that
 \begin{equation}
 \label{E2.26nnn}
\begin{array}{l}
  \di\sum_{k\in\N} b(j_0(k))k^s\| \na u^{k/2}\|_2^2\leq \|\na u^{1/2}\|_2^2\\[3ex]
\qquad\di +\sum_{n=1}^{n_0-1} (3n)^s\|\na u_{n,n+1}\|_2^2
\sum_{k=1}^n
\big(b(j_0(2k))+b(j_0(2k+1))\big)\\[3ex]
\qquad  +\di\sum_{n\geq n_0} (3n)^s\|\na u_{n,n+1}\|_2^2
 \sum_{k=1}^n \big(b(j_0(2k))+b(j_0(2k+1))\big)\\[3ex]
 \quad\di \leq \|\na u^{1/2}\|_2^2+(c(n_0,s)-1)\|\na u_{n_0}\|_2^2+2\cdot 3^{s+1}\sum_{n\geq n_0}n^{s+1}\|\na u_{n,n+1}\|_2^2
\\[3ex]
 \quad\di \leq c(n_0,s)\|\na u\|_2^2+2\cdot 3^{s+1}\sum_{n\geq n_0} n^{s+1}\|\na u_{n,n+1}\|_2^2,
   \end{array}
  \end{equation}
  where $c(n_0,s)$ is a generic constant depending only on $n_0$ and $s$.

Therefore, we get from \eq{2.5n} -- \eq{2.26nnn} that
 \begin{equation}
  \label{E4.13n}
  \begin{array}{l}
 \di\frac{d}{2dt}\di\sum_{k\in\N}k^s\|u^{k}\|_2^2+\nu \sum_{k\in\N}k^s\|\na u^{k}\|_2^2  \\[3ex]
 \leq 4C_1\|u(t)\|_{X_1}\big(c(n_0,s)\|\na u(t)\|_2^2
   +2\cdot 3^{s+1}\di \sum_{n\geq n_0} n^{s+1}\|\na u_{n,n+1}(t)\|_2^2\big)\ek
\leq C_2\|u(t)\|_{X_1}\big(\|\na u(t)\|_2^2
   +\di \sum_{n\in\N} n^{s+1}\|\na
   u_{n,n+1}(t)\|_2^2\big),\quad\text{in }
   (0,T),
   \end{array}
 \end{equation}
where $C_2$ is a generic constant depending only on $s$ and $d$.

By integrating \eq{4.13n} from $0$ to $t\in (0,T)$  we have
 \begin{equation}
 \label{E4.12n}
 \begin{array}{l}
\di\ha\sum_{k\in\N}k^s\|u^{k}(t)\|_2^2+\nu \int_0^t
\sum_{k\in\N}k^s\|\na u^{k}(\tau)\|_2^2\,d\tau
\leq \ha\di\sum_{k\in\N}k^s\|[u_0]^{k}\|_2^2\\[2ex]
  \quad +C_2\di\int_0^t \|u(\tau)\|_{X_1}\big(\|\na u(\tau)\|_2^2
   +\di \sum_{n\in\N} n^{s+1}\|\na
   u_{n,n+1}(\tau)\|_2^2\big)\,d\tau,\,\forall t\in (0,T).
  \end{array}
  \end{equation}
Note that
 \begin{equation}
 \label{E2.25nn}
\begin{array}{l} \di \sum_{k\in \N}k^s\|u^{k}\|_2^2
=\di\sum_{n\in\N}\big(\sum_{k=1}^n k^s\big)\|u_{n,n+1}\|_2^2,\ek
 \di \sum_{k\in\N}\|\na u^{k}\|_2^2 =\di\sum_{n\in\N}\big(\sum_{k=1}^n
 k^s\big)\|u_{n,n+1}\|_2^2,
 \end{array}
 \end{equation}
 where
 \begin{equation}
 \label{E2.25nnn}
\frac{n^{s+1}}{s+1}\leq \sum_{k=1}^n k^s\leq
\frac{(n+1)^{s+1}}{s+1}.
 \end{equation}
By \eq{4.12n}--\eq{2.25nnn} we have
 \begin{equation}
 \label{E4.12nn}
 \begin{array}{l}
\di\ha\sum_{n\in\N}n^{s+1}\|u_{n,n+1}(t)\|_2^2+\nu \int_0^t
\sum_{n\in\N}n^{s+1}\|\na u_{n,n+1}(\tau)\|_2^2\,d\tau\\[2ex]
\quad \leq \di\ha\sum_{n\in\N}(n+1)^{s+1}\|[u_0]_{n,n+1}\|_2^2\\[2ex]
  \; +C_2(s+1)\di\int_0^t \|u(\tau)\|_{X_1}\big(\|\na u(\tau)\|_2^2
   +\di \sum_{n\in\N} n^{s+1}\|\na u_{n,n+1}(\tau)\|_2^2\big)\,d\tau,
   \, \forall t\in (0,T).
  \end{array}
  \end{equation}
 Note that,
in view of $|\xi|^{s+1}-n^{s+1}\leq (s+1)|\xi|^s$ for $n<|\xi|\leq
n+1$,
\begin{equation}
 \label{E3.4nn}
\ha\|u^1(t)\|_{\dot{H}^{(s+1)/2}}^2-\ha\sum_{n\in\N}n^{s+1}\|u_{n,n+1}(t)\|_2^2
        \leq \frac{s+1}{2}\|u^1(t)\|_{\dot{H}^{s/2}}^2,
        \end{equation}
\begin{equation}
 \label{E3.4nnn}
 \begin{array}{l}
\nu\di\int_0^t\Big(\|\na
u^1(\tau)\|_{\dot{H}^{(s+1)/2}}^2-\sum_{n\in\N}n^{s+1}\|\na
u_{n,n+1}(\tau)\|_2^2\Big)\,d\tau\\[2ex]
\qquad \qquad\qquad\qquad\leq (s+1)\nu\di\int_0^t\|\na
u^1(\tau)\|_{\dot{H}^{s/2}}^2\,d\tau.
 \end{array}
 \end{equation}
Thus, summing up \eq{4.12nn}, \eq{3.4nn},  \eq{3.4nnn} and taking
into account
 $$\begin{array}{l}
\di\sum_{n\in\N}(n+1)^{s+1}\|[u_0]_{n,n+1}\|_2^2\leq
2^{s+1}\|u_0\|_{\dot{H}^{(s+1)/2}}^2,\ek
  \di\sum_{n\in\N}n^{s+1}\|\na u_{n,n+1}\|_2^2\leq\|u^1\|_{\dot{H}^{(s+3)/2}}^2,
  \end{array}
$$
 we obtain
 \begin{equation}
 \label{E4.20}
 \begin{array}{l}
\di\ha\|u^1(t)\|_{\dot{H}^{(s+1)/2}}^2+\nu \int_0^{t}
 \|\na u^1(\tau)\|_{\dot{H}^{(s+1)/2}}^2\,d\tau\\[2ex]
\quad \leq \di(s+1)\big(\ha\|u^1(t)\|_{\dot{H}^{s/2}}^2+\nu
\int_0^{t}
 \|\na u^1(\tau)\|_{\dot{H}^{s/2}}^2\,d\tau\big)
+2^{s}\|u_0\|_{\dot{H}^{(s+1)/2}}^2\\[2ex]
  \qquad +C_3\di\int_0^t \|u(\tau)\|_{X_1}\big(\|\na u(\tau)\|_2^2
   + \|\na u^1(\tau)\|_{\dot{H}^{(s+1)/2}}^2\big)\,d\tau,
   \; \forall t\in (0,T),
  \end{array}
  \end{equation}
  where $C_3:=C_2(s+1)$.
Since  the inequality
 $$\ha\|u_1(t)\|_{\dot{H}^{(s+1)/2}}^2+\nu
\int_0^{t} \|\na u_1(\tau)\|_{\dot{H}^{(s+1)/2}}^2\,d\tau\leq
\ha\|u_0\|_2^2,\,\forall t\in (0,T),$$ holds thanks to
$\|v_1\|_{\dot{H}^{(s+1)/2}}\leq \|v_1\|_2\leq \|v\|_2,\forall v\in
H^{(s+1)/2}(\R^d),$ and the energy inequality \eq{EI} in $(0,T)$, we
get by \eq{4.20} that
\begin{equation}
 \label{E4.21}
 \begin{array}{l}
\di\ha\|u(t)\|_{\dot{H}^{(s+1)/2}}^2+\nu \int_0^{t}
 \|\na u(\tau)\|_{\dot{H}^{(s+1)/2}}^2\,d\tau\\[2ex]
  \leq \di\ha\|u_1(t)\|_{\dot{H}^{(s+1)/2}}^2+\nu \int_0^{t}
 \|\na u_1(\tau)\|_{\dot{H}^{(s+1)/2}}^2\,d\tau\ek
  \quad+\di(s+1)\big(\ha\|u^1(t)\|_{\dot{H}^{s/2}}^2+\nu
\int_0^{t}
 \|\na u^1(\tau)\|_{\dot{H}^{s/2}}^2\,d\tau\big)
+2^s\|u_0\|_{\dot{H}^{(s+1)/2}}^2\\[2ex]
  \qquad +C_3\di\int_0^t \|u(\tau)\|_{X_1}\big(\|\na u(\tau)\|_2^2
   + \|\na u^1(\tau)\|_{\dot{H}^{(s+1)/2}}^2\big)\,d\tau\ek
     \leq \di\ha\|u_0\|_{2}^2+2^s\|u_0\|_{\dot{H}^{(s+1)/2}}^2+\di(s+1)\big(\ha\|u(t)\|_{\dot{H}^{s/2}}^2+\nu
\int_0^{t}
 \|\na u(\tau)\|_{\dot{H}^{s/2}}^2\,d\tau\big)
\\[2ex]
  \qquad +C_3\di\int_0^t \|u(\tau)\|_{X_1}\big(\|\na u(\tau)\|_2^2
   + \|\na u(\tau)\|_{\dot{H}^{(s+1)/2}}^2\big)\,d\tau,\;\forall t\in (0,T).
  \end{array}
  \end{equation}
 Using $$(s+1)\|v\|_{\dot{H}^{s/2}}^2\leq c(s)\|v\|_2^2+\ha\|v\|_{\dot{H}^{(s+1)/2}}^2,
  \forall v\in H^{(s+1)/2}(\R^d),$$
  due to complex interpolation
$\dot{H}^{s/2}=[L^2, \dot{H}^{(s+1)/2}]_{s/(s+1)}$, we get from
\eq{4.21}
 that
\begin{equation}
 \label{E4.23}
 \begin{array}{l}
\di\frac{1}{4}\|u(t)\|_{\dot{H}^{(s+1)/2}}^2+\frac{\nu}{2}\int_0^{t}
 \|\na u(\tau)\|_{\dot{H}^{(s+1)/2}}^2\,d\tau\\[2ex]
     \leq \di\ha\|u_0\|_{2}^2+2^s\|u_0\|_{\dot{H}^{(s+1)/2}}^2+\di c(s)\big(\ha\|u(t)\|_{2}^2+\nu
\int_0^{t} \|\na u(\tau)\|_{2}^2\,d\tau\big)
\\[2ex]
  \qquad +C_3\di\int_0^t \|u(\tau)\|_{X_1}\big(\|\na u(\tau)\|_2^2
   + \|\na u(\tau)\|_{\dot{H}^{(s+1)/2}}^2\big)\,d\tau\\[2ex]
      \leq \di\ha(1+c(s))\|u_0\|_{2}^2+2^s\|u_0\|_{\dot{H}^{(s+1)/2}}^2
\\[2ex]
  \qquad +C_3\di\int_0^t \|u(\tau)\|_{X_1}\big(\|\na u(\tau)\|_2^2
   + \|\na u(\tau)\|_{\dot{H}^{(s+1)/2}}^2\big)\,d\tau,\;\forall t\in (0,T).
  \end{array}
  \end{equation}

Now, put $\ve=\frac{\nu}{4C_3}$ and fix $t\in (0,T)$. By Lemma
\ref{L4.2n} (ii) in view of $u\in C([0,t],
  \dot{H}^{-1+d/2}(\R^d))$ there exists
$l_0(t)>0$ such that
 \begin{equation}
 \label{E4.16}
\begin{array}{l}
  \di\|\la u(\tau, \la\cdot)\|_{X_1}\leq \ve,\quad  \forall \tau\in [0,t],\;\forall
\la=2^{2^{2^l}} (l\geq l_0(t)).
 \end{array}
 \end{equation}
Put $\la=2^{2^{2^l}}$ with fixed $l\geq l_0(t)$  and observe
  the re-scaled function
 \begin{equation}
 \label{E4.19}
 (u)_\la(\tau,x):=\la u(\la^2 \tau, \la x), \,\tau\in
 (0,\la^{-2}t].
  \end{equation}
Then, \eq{4.16} implies that
 \begin{equation}
 \label{E4.17}
\begin{array}{l}
 \di \|(u)_\la(\tau)\|_{X_1}\leq \ve,\quad  \forall \tau\in [0,\la^{-2}t].
 \end{array}
 \end{equation}

Since $(u)_\la$ is the Leray-Hopf weak solution to \eq{1.1} with
initial value $(u_0)_\la$ and the first blow-up epoch $\la^{-2}T$,
we have by  \eq{4.23} that
 \begin{equation}
 \label{E4.15}
 \begin{array}{l}
\di\frac{1}{4}\|(u)_\la(\la^{-2}t)\|_{\dot{H}^{(s+1)/2}}^2+\frac{\nu}{2}
\int_0^{\la^{-2}t}
 \|\na (u)_\la(\tau)\|_{\dot{H}^{(s+1)/2}}^2\,d\tau\\[2ex]
\leq \di\ha(1+c(s))\|(u_0)_\la\|_2^2+2^s\|(u_0)_\la\|_{\dot{H}^{(s+1)/2}}^2\\[2ex]
  \quad +C_3\di\int_0^{\la^{-2}t} \|(u)_\la(\tau)\|_{X_1}\big(\|\na (u)_\la(\tau)\|_2^2
   +\|\na (u)_\la(\tau)\|_{\dot{H}^{(s+1)/2}}^2\big)\,d\tau,
  \end{array}
  \end{equation}
which yields by \eq{4.17} and energy inequality for $(u)_\la$ in
$(0,\la^{-2}T)$ that
 \begin{equation}
 \label{E4.18n}
 \begin{array}{l}
\di\frac{1}{4}\|(u)_\la(\la^{-2}t)\|_{\dot{H}^{(s+1)/2}}^2+\frac{\nu}{4}\int_0^{\la^{-2}t}
 \|\na (u)_\la(\tau)\|_{\dot{H}^{(s+1)/2}}^2\,d\tau
\\[2ex]
  \qquad \di\leq \ha(1+c(s))\|(u_0)_\la\|_2^2+2^s\|(u_0)_\la\|_{\dot{H}^{(s+1)/2}}^2+\frac{\nu}{4} \int_0^{\la^{-2}t} \|\na
  (u)_\la(\tau)\|_2^2\,d\tau\\[2ex]
  \qquad \di\leq (1+c(s))\|(u_0)_\la\|_2^2+2^s\|(u_0)_\la\|_{\dot{H}^{(s+1)/2}}^2.
  \end{array}
  \end{equation}
Thanks to the scaling-variant property
$$\begin{array}{l}
\|(u)_\la(\la^{-2}t)\|_{\dot{H}^{(s+1)/2}}^2=\la^{3+s-d}
\|u(t)\|_{\dot{H}^{(s+1)/2}}^2, \\[2ex]
 \int_0^{\la^{-2}t}\|\na (u)_\la(\tau)\|_{\dot{H}^{(s+1)/2}}^2\,d\tau=\la^{3+s-d}\int_0^{t}
 \|\na u(\tau)\|_{\dot{H}^{(s+1)/2}}^2\,d\tau,\\[2ex]
 \|(u_0)_\la\|_{\dot{H}^{(s+1)/2}}^2=\la^{3+s-d}\|u_0\|_{\dot{H}^{(s+1)/2}}^2,
  \|(u_0)_\la\|_2^2=\la^{2-d}\|u_0\|_2^2,
 \end{array}$$
 we get by \eq{4.18n} that
$$ \begin{array}{l} \di\|u(t)\|_{\dot{H}^{(s+1)/2}}^2+\nu\int_0^{t}
 \|\na u(\tau)\|_{\dot{H}^{(s+1)/2}}^2\,d\tau\ek
  \qquad\qquad\leq
  4(1+c(s))\la^{-s-1}\|u_0\|_2^2+2^{s+2}\|u_0\|_{\dot{H}^{(s+1)/2}}^2.
  \end{array}
  $$
Consequently, in view of $\la\geq 4$,  we have
$$
\begin{array}{l} \di\|u(t)\|_{\dot{H}^{(s+1)/2}}^2+\nu\int_0^{t}  \|\na
u(\tau)\|_{\dot{H}^{(s+1)/2}}^2\,d\tau \leq
(1+c(s))\|u_0\|_2^2+2^{s+2}\|u_0\|_{\dot{H}^{(s+1)/2}}^2,
  \end{array}
  $$
which together with the energy inequality \eq{EI} in $(0,T)$ yields
that
 \begin{equation}
 \label{E4.24}
\begin{array}{l} \di\|u(t)\|_{H^{(s+1)/2}}^2+\nu\int_0^{t}  \|\na
u(\tau)\|_{H^{(s+1)/2}}^2\,d\tau \leq C(s)\|u_0\|_{H^{(s+1)/2}}^2.
  \end{array}
  \end{equation}
Note here that $C(s)$ is irrespective of $d$ and $t$, of course.

Thus we have obtained \eq{4.24} for any fixed $t\in (0,T)$.
 Since $\frac{s+1}{2}\geq -1+\frac{d}{2}$ for $s\geq d-3$, we can conclude from \eq{4.24}  that
$T=\infty$. Moreover, \eq{4.24} with $s\geq d-3$ implies \eq{1.2n}
with $s\geq -1+d/2$.

The proof of the theorem is complete.

\hfill\qed
\begin{rem}
 {\rm
 The above argument essentially requires  that the initial value
$u_0$ should be taken, at least, in $H^{-1+d/2}(\R^d)$ so that a
local smooth solution can exist and the uniform smallness condition
\eq{4.17} for the re-scaled solution can be satisfied.
 }
\end{rem}

\section{Application to global solvability for incompressible Euler equations}

\par\bigskip

 The above regularity result for the  Navier-Stokes system
 can be directly applied to showing the global existence of solutions to
 the Euler equations of ideal incompressible fluids:
\begin{equation}
\label{E1.4}
\begin{array}{rl}
     u_t+ (u\cdot\na)u+\na p  = 0 \,\, &\text{in }(0,\infty)\ti\R^d,\ek
      \div u = 0 \,\, &\text{in }(0,\infty)\ti\R^d,\ek
u(0,x)=u_0\, &\text{in }\R^d\; (d\geq 3).
\end{array}
\end{equation}
\par\medskip
Let us denote by $a\otimes b:=(a_ib_j)$ for vectors $a,b\in \R^d$
and $A:B=\sum_{i,j}a_{ij}b_{ij}$ for matrices $A=(a_{ij}),
B=(b_{ij})\in \R^{d\times d}$. A $d$-dimensional divergence-free
vector field $u\in L_{\rm loc}^2([0,\infty)\ti \R^d)$ is called a
weak solution to \eq{1.4} if it satisfies the equations in a weak
sense, i.e.,
 \begin{equation}
 \label{E4.5}
\begin{array}{l}
\int_0^\infty \int_{\R^d} (u\cdot \vp_t+(u\otimes u): \na\vp)dxdt
 +\int_{\R^d} u_0(x)\cdot\vp(0,x)dx=0,\ek\hfill
\forall \vp\in C_0^\infty([0,\infty)\ti\R^d)\; (\div\vp=0).
 \end{array}
 \end{equation}
\begin{rem}
 \label{R4.1}
  {\rm
  (i) If  $u$ is a weak solution
 to \eq{1.4}, then by De Rham's theorem there is a distribution $p$ such
that $\{u ,p\}$ satisfies the equations of \eq{1.4} in the
distribution sense.

(ii) Suppose a weak solution $u$ to \eq{1.4} belongs to
$L^\infty(0,\infty;H^s(\R^d))$ for $s>1+\frac{d}{2}$.  It then
follows by the Sobolev embedding $H^{s-1}(\R^d)\hra L^\infty(\R^d)$
that $(u\cdot\na)u\in L^\infty(0,\infty; L^2(\R^d))$ and   $$-\Da
p=\div (u\cdot\na)u\in L^\infty(0,\infty;H^{-1}(\R^d))$$
 holds for the
corresponding associated pressure $p$. Consequently,  $\na p$ and
hence $u_t$ belong to $L^\infty(0,\infty; L^2(\R^d))$.
  Therefore, $u$ becomes a strong solution to \eq{1.4}.
 }
\end{rem}
\par\medskip
For \eq{1.4} in 3D setting, nonuniqueness of weak solutions  for
some initial values is proved (e.g. in \cite{DeLeSz10},\cite{Is17})
and local existence of a unique strong solution for $u_0\in
H^s(\R^3), s>5/2,$ is known (cf. \cite{BaTi07}). However, there
seems no general establishment asserting global existence of a
solution under a certain regularity assumption on the initial
values.
\par\medskip
As a corollary of  Theorem \ref{T1.1}, we can obtain the following
global solvability for \eq{1.4}.
\begin{theo}
\label{T1.3} {\rm Let $u_0\in H^{s}(\R^d)$, $d\geq 3$, $s\geq
-1+d/2$, $\div u_0=0$.
 Then the Cauchy problem \eq{1.4} has a global weak solution such
 that
$$u\in
L^\infty(0,\infty;H^{s}(\R^d))
$$
 and
 \begin{equation}
 \label{E1.7}
\|u\|_{L^\infty(0,\infty;H^{s}(\R^d))}\leq C\|u_0\|_{H^s(\R^d)}
\end{equation}
with some constant $C>0$. This solution is unique if
  $s>1+d/2$.}
\end{theo}
\par\medskip
 {\bf Proof:}
The proof is based on the argument of Navier-Stokes regularization
of the problem \eq{1.4} and vanishing viscosity.
 Observe the
Navier-Stokes problem \eq{1.1} with the same initial value $u_0$.

By Theorem \ref{T1.1}, the Navier-Stokes problem \eq{1.1} has a
unique solution $u_\nu \in L^\infty(0,\infty;H^{s}(\R^d))$
satisfying
 \begin{equation}
 \label{E4.6}
\begin{array}{l}
\int_0^\infty \int_{\R^d} (u_\nu\cdot \vp_t-\nu\na u_\nu
\cdot\na\vp+(u_\nu\otimes u_\nu): \na\vp)dxdt
 +\int_{\R^d} u_0(x)\cdot\vp(0,x)dx=0,\ek\hfill
\forall \vp\in C_0^\infty([0,\infty)\ti\R^d)\; (\div\vp=0),
 \end{array}
 \end{equation}
and
\begin{equation}
\label{E4.2} \|u_\nu\|_{L^\infty(0,\infty;H^{s}(\R^d))}\leq
C(s)\|u_0\|_{H^s(\R^d)},
 \end{equation}
 where $C(s)$  is independent of $\nu$.

Let $u\in L^\infty(0,\infty;H^{s}(\R^d))$ be the  weak-$\ast$ limit
of a subsequence of $\{u_\nu\}_{\nu>0}$ (say, $\{u_\nu\}_{\nu>0}$
again).  Since
$$\begin{array}{l}\nu\big|\int_0^\infty \int_{\R^d} \na u_\nu
\cdot\na\vp\,dxdt\big|\leq \nu
\|u_\nu\|_{L^\infty(0,\infty;H^{s}(\R^d))}\|\Da\vp\|_{L^1(0,\infty;H^{-s}(\R^d))}\ek
\qquad \leq \nu
C\|u_0\|_{H^{s}(\R^d)}\|\Da\vp\|_{L^1(0,\infty;H^{-s}(\R^d))}\ra
0\quad (\nu\ra 0)
\end{array}$$
 for all $\vp\in C_0^\infty([0,\infty)\ti\R^d)$,
 the integral identity \eq{4.5} will be satisfied by $u$
  provided
 \begin{equation}
 \label{E4.8}
\begin{array}{l}
\int_0^\infty \int_{\R^d}(u_\nu\otimes u_\nu): \na\vp dxdt\ra
\int_0^\infty \int_{\R^d}(u\otimes u): \na\vp dxdt\quad (\nu\ra
0),\ek \hfill \forall \vp\in C_0^\infty([0,\infty)\ti\R^d)\;
(\div\vp=0).\end{array}\end{equation}

Let a distribution $p_\nu$  be an associated pressure for \eq{1.1}
corresponding to $u_\nu$, that is, the pair $\{u_\nu,p_\nu\}$
satisfies the momentum equation of \eq{1.1} in the distribution
sense. Since
 \begin{equation}
 \label{E4.7}
H^{-1+d/2}(\R^d)\cdot H^{-1+d/2}(\R^d)\hra H^{-1/2}(\R^d), d\geq 3,
 \end{equation}
   it follows by \eq{4.2} that
 $\{u_\nu\otimes u_\nu\}_{\nu>0}$ is uniformly bounded in $L^\infty(0,\infty;H^{-1/2}(\R^d))$.
On the other hand, since the associated pressure $p_\nu$  satisfies
 $$-\Da p_\nu=\div \div (u_\nu\otimes u_\nu),$$
 we have, in view of \eq{4.7},
  $$\begin{array}{l}\|p_\nu\|_{L^\infty(0,\infty;H^{-1/2}(\R^d))}\leq C \|u_\nu\otimes
 u_\nu\|_{L^\infty(0,\infty;H^{-1/2}(\R^d))}\ek
 \qquad \leq C\|u_\nu\|_{L^\infty(0,\infty;H^{-1+d/2}(\R^d))}^2\leq
 C\|u_0\|_{H^{-1+d/2}(\R^d)}^2\end{array}$$
with $C>0$ independent of $\nu>0$. Therefore, in view of $$u_{\nu
t}=\Da u_\nu - \div (u_\nu\otimes u_\nu)-\na p_\nu,$$ it follows
that $\{u_{\nu t}\}_{\nu>0}$ is uniformly bounded in
$L^\infty(0,\infty; H^{-3/2}(\R^d))$.

Since  $\{u_{\nu}\}_{\nu>0}$ and $\{u_{\nu t}\}_{\nu>0}$ are
uniformly bounded in $L^\infty(0,\infty;H^{s}(\R^d))$ and in
$L^\infty(0,\infty; H^{-3/2}(\R^d))$, respectively, the restriction
of both sequences on $(0,T)\ti G$  for any finite $T>0$ and smooth
bounded domain $G\subset \R^d$ are uniformly bounded in
$L^2(0,T;H^{s}(G))$ and in $L^2(0,T; H^{-3/2}(G))$, respectively,
Therefore, it follows by Aubin's lemma on compactness that, possibly
with a new subsequence,
$$u_\nu \ra u \quad\text{in } L^2(0,T; L^2(G)) \quad\text{as }\nu\ra 0 \quad\text{(strongly)}.$$
 Consequently, \eq{4.8} holds and $u$ is a weak solution satisfying
\eq{4.5} and \eq{1.7}.

Let us prove the uniqueness of the solution in
$L^\infty(0,\infty;H^s(\R^d))$ for $s>1+\frac{d}{2}$. Let
$u_1,u_2\in L^\infty(0,\infty;H^s(\R^d))$, $s>1+\frac{d}{2}$, be two
solutions to \eq{1.4}, and let $u\equiv u_1-u_2$ and $p\equiv
p_1-p_2$ with corresponding associated pressure $p_1$, $p_2$. Then,
$$u_t+ (u\cdot\na) u_1+ (u_2\cdot\na) u+\na p  = 0\quad\text{in }(0,\infty).$$
Hence, in view of Remark \ref{R4.1} (ii), we have
$$(u_t,u)_{L^2}+ \big((u\cdot \na)u_1, u\big)_{L^2}= 0\quad\text{in }(0,\infty),$$
which yields
$$\frac{d}{2dt}\|u(t)\|_{2}^2\leq \|\na u_1(t)\|_\infty \|u(t)\|_{2}^2
 \leq c\|u_1(t)\|_{H^s(\R^d)} \|u(t)\|_{2}^2,\; \forall t\in (0,\infty).$$
Thus, by  Gronwall's inequality we have $u(t)\equiv 0, \forall t\in
(0,\infty)$.

The proof of the theorem is complete. \hfill\qed


\bigskip

\bigskip\noindent
\textbf{Conflict of interest:} The author declares that there is no
potential conflict of interest with other people or organization
associated with this manuscript.

\end{document}